# DYNAMIC EXPONENTIAL UTILITY INDIFFERENCE VALUATION

By Michael Mania[1] and Martin Schweizer

*A. Razmadze Mathematical Institute and ETH Zürich*

We study the dynamics of the exponential utility indifference value process $C(B; \alpha)$ for a contingent claim $B$ in a semimartingale model with a general continuous filtration. We prove that $C(B; \alpha)$ is (the first component of) the unique solution of a backward stochastic differential equation with a quadratic generator and obtain *BMO* estimates for the components of this solution. This allows us to prove several new results about $C_t(B; \alpha)$. We obtain continuity in $B$ and local Lipschitz-continuity in the risk aversion $\alpha$, uniformly in $t$, and we extend earlier results on the asymptotic behavior as $\alpha \searrow 0$ or $\alpha \nearrow \infty$ to our general setting. Moreover, we also prove convergence of the corresponding hedging strategies.

**0. Introduction.** One of the important problems in mathematical finance is the valuation of contingent claims in incomplete financial markets. In mathematical terms, this can be formulated as follows. We have a semimartingale $S$ modeling the discounted prices of the available assets and a random variable $B$ describing the payoff of a financial instrument at a given time $T$. The gains from a trading strategy $\vartheta$ with initial capital $x$ are described by the stochastic integral $x + \int \vartheta \, dS = x + G(\vartheta)$. If $B$ admits a representation as $B = x + G_T(\vartheta)$ for some pair $(x, \vartheta)$, the claim $B$ is called attainable, and its value at any time $t \leq T$ must equal $x + G_{0,t}(\vartheta)$ due to absence-of-arbitrage considerations. Incompleteness means that there are some nonattainable $B$, and the question is how to value those.

In this paper we use the utility indifference approach to this problem. For a given utility function $U$ and an initial capital $x_t$ at time $t$, we define the value $C_t(x_t, B)$ implicitly by the requirement that

$$\operatorname*{ess\,sup}_{\vartheta} E[U(x_t + G_{t,T}(\vartheta)) | \mathcal{F}_t] = \operatorname*{ess\,sup}_{\vartheta} E[U(x_t + C_t(x_t, B) + G_{t,T}(\vartheta) - B) | \mathcal{F}_t].$$

Received September 2004.
[1]Supported in part by the Grant Lagrange-CRT Foundation-ISI Foundation.
*AMS 2000 subject classifications.* 91B28, 60H10, 91B16, 60G48.
*Key words and phrases.* Indifference value, exponential utility, dynamic valuation, BSDE, semimartingale backward equation, *BMO*-martingales, incomplete markets, minimal entropy martingale measure.

In terms of expected utility, we are thus indifferent between selling or not selling the claim $B$ for $C_t(x_t, B)$, provided that we combine each of those alternatives with optimal trading. Our goal is to study the dynamic behavior of the process $C = C(B; \alpha)$ resulting from the exponential utility function $U(x) = -e^{-\alpha x}$ with risk aversion $\alpha \in (0, \infty)$.

The existing literature on exponential utility indifference valuation can be roughly divided in two groups. A larger set of papers studies static questions; they examine $C_0(B; \alpha)$, the time 0 value, in models of varying generality. A good recent overview with many references is given by Becherer [3]; [17] contains a slightly different approach and additional references. The second set of papers studies $C(B; \alpha)$ as a process; this is done by Rouge and El Karoui [33] in a Brownian filtration, by Musiela and Zariphopoulou [29] or Young [35], among several others, in a Markovian diffusion setting or by Musiela and Zariphopoulou [30] in a binomial model. In the present paper, we work in a general continuous filtration and obtain several new results on the dynamic properties of the process $C(B; \alpha)$ and its asymptotic behavior as the risk aversion $\alpha$ goes to 0 or to $\infty$. In particular, we provide convergence results for hedging strategies.

The paper is structured as follows. Section 1 lays out the model and provides some auxiliary results mostly known from the literature. We use these to represent $C(B; \alpha)$ as the dynamic value process of a standard utility maximization problem with a random endowment and formulated under a suitable measure $Q^E$. This allows us, in Section 2, to extend the static properties known for $C_0(B; \alpha)$ very easily to any $C_t(B; \alpha)$. Moreover, we easily obtain the existence of an optimal strategy for this stochastic control problem. Section 4 shows that $C(B; \alpha)$ is the unique solution of a backward stochastic differential equation (BSDE) with a quadratic generator. In contrast to a similar result by Rouge and El Karoui [33], our derivation directly uses the martingale optimality principle and the existence of an optimal strategy. Section 3 prepares for these results by providing a comparison theorem for a more general class of BSDEs driven by a martingale in a continuous filtration and having quadratic generators. The key step here is Proposition 7, which shows that the martingale part of any bounded solution of a BSDE with a generator satisfying a quadratic growth condition belongs to $BMO$. This underlines the importance of $BMO$-martingales when dealing with BSDEs with quadratic generators. For the particular generator corresponding to the BSDE for $C(B; \alpha)$, we also obtain estimates on the $BMO$ norms of the components of the solution.

Section 5 exhibits additional properties of the valuation $C(B; \alpha)$. We obtain time-consistency, continuity in $B$ and local Lipschitz-continuity in $\alpha$, both of the latter uniformly in $t$. Finally, Section 6 studies the asymptotics of $C(B; \alpha)$ as $\alpha$ goes to 0 or to $\infty$. For $\alpha \searrow 0$, we prove generally that $C_t(B; \alpha)$ decreases to $E_{Q^E}[B|\mathcal{F}_t]$ at a rate of $\alpha$, uniformly in $t$, where $Q^E$ is the minimal entropy martingale measure for $S$; this is a simple extension of a result due to Stricker [34]. With the help of our BSDE description, we are, moreover, able to prove the novel result that the corresponding hedging strategies $\psi(\alpha)$ converge to the strategy $\psi^E$ which is risk-minimizing under $Q^E$ in the sense of Föllmer and Sondermann [15]. For $\alpha \nearrow \infty$, $C_t(B;\alpha)$ increases generally to the superreplication price $C_t^*(B) = \operatorname{ess\,sup}_Q E_Q[B|\mathcal{F}_t]$, uniformly in $t$; this generalizes a result due to Rouge and El Karoui [33] for the case of a Brownian filtration. In addition, again using the BSDE, we also prove the convergence of the corresponding hedging strategies $\psi(\alpha)$ to the superreplication strategy $\psi^*$ from the optional decomposition of $C^*(B)$.

**1. Basic concepts and preliminary results.** In this section we introduce the notion of the utility indifference value process for a contingent claim and recall some basic facts for the case of an exponential utility function.

We start with a probability space $(\Omega, \mathcal{F}, P)$, a time horizon $T \in (0, \infty]$ and a filtration $\mathbb{F} = (\mathcal{F}_t)_{0 \leq t \leq T}$ satisfying the usual conditions of right-continuity and completeness. Hence, we can and do choose RCLL versions for all semimartingales. Fix an $\mathbb{R}^d$-valued semimartingale $S = (S_t)_{0 \leq t \leq T}$ and think of this as the *discounted price process* for $d$ risky assets in a financial market containing also a riskless asset with discounted price constant at 1. A *self-financing trading strategy* is determined by its initial capital $x \in \mathbb{R}$ and the numbers $\vartheta_t^i$ of shares of asset $i$, $i = 1, \ldots, d$, held at time $t \in [0, T]$. Formally, $\vartheta$ is in the space $L(S)$ of $\mathbb{F}$-predictable $S$-integrable $\mathbb{R}^d$-valued processes so that the (real-valued) stochastic integrals $G_{t,u}(\vartheta) := \int_t^u \vartheta_s \, dS_s$ are well defined. They describe the gains or losses from trading according to $\vartheta$ between $t$ and $u > t$. The wealth at time $t$ of a strategy $(x, \vartheta)$ is $x + G_{0,t}(\vartheta)$ and we denote by $G(\vartheta)$ the running stochastic integral process $G_{0,\cdot}(\vartheta)$. Arbitrage opportunities will be excluded below via the choice of a suitable space $\Theta$ of "permitted" trading strategies $\vartheta$.

Now let $U : \mathbb{R} \to \mathbb{R}$ be a utility function and $B \in L^0(\mathcal{F}_T)$ a *contingent claim*, that is, a random payoff at time $T$ described by the $\mathcal{F}_T$-measurable random variable $B$. In order to assign to $B$ at some date $t \in [0, T]$ a (subjective) value based on the utility function $U$, we first fix an $\mathcal{F}_t$-measurable random variable $x_t$. Then we define

$$V_t^B(x_t) := \operatorname*{ess\,sup}_{\vartheta \in \Theta} E[U(x_t + G_{t,T}(\vartheta) - B)|\mathcal{F}_t],$$

the maximal conditional expected utility we can achieve by starting at time $t$ with initial capital $x_t$, using some strategy $\vartheta \in \Theta$ on $(t, T]$ and paying out $B$ at time $T$. The *utility indifference value* $C_t(x_t, B)$ at time $t$ for $B$ with respect to $U$ and $x_t$ is implicitly defined by

$$(1.1) \qquad V_t^0(x_t) = V_t^B(x_t + C_t(x_t, B)).$$

This says that starting with $x_t$, one has the same maximal utility from solely trading on $(t, T]$ as from selling $B$ at time $t$ for $C_t(x_t, B)$, again trading and then paying out $B$ at the final date $T$.

REMARK. Variants of the above notion of utility indifference value have been known and used for a long time. Its first appearance in a form that also accounts for the presence of a financial market is usually attributed to Hodges and Neuberger [20]. The resulting valuation has been studied extensively in recent years and we shall provide some more references when giving more specific results. One good starting point with a long literature list is [3]. However, most papers only define this value for $t = 0$ and with $\mathcal{F}_0$ trivial and thus obtain one mapping from (random payoffs $B$ in) $L^0(\mathcal{F}_T)$ to $\mathbb{R}$. Exceptions are papers set in Markovian frameworks where the stochastic processes $V^B(x)$ and $C(x, B)$ can be represented via functions of the state variables; see, for instance, [29] or [35] for recent papers with more references to earlier work. There is also some literature on dynamic versions of this valuation and their properties; see, notably, [33] or [3, 4]. But in contrast to our approach, these authors use the definition (1.1) only for $t = 0$ and another one for $t \in (0, T]$, and they do not argue (the fact) that their definitions are equivalent to (1.1) for all $t$.

To pass from the above formal definitions to rigorous results, we now choose one particular $U$ and a corresponding $\Theta$. Throughout the rest of this paper, we work with the *exponential utility* function

$$U(x) = -\exp(-\alpha x)$$

with risk aversion $\alpha \in (0, \infty)$. We assume that

$$S \text{ is locally bounded,}$$

denote by $\mathbb{P}_e := \{Q \approx P | S \text{ is a local } Q\text{-martingale}\}$ the set of all equivalent local martingale measures for $S$ and assume that

$$\mathbb{P}_{e,f} := \{Q \in \mathbb{P}_e | H(Q|P) < \infty\} \neq \varnothing.$$

Finally, we define the space of our trading strategies as

$$\Theta := \{\vartheta \in L(S) | G(\vartheta) \text{ is a } Q\text{-martingale for all } Q \in \mathbb{P}_{e,f}\}.$$

For future use, we introduce the terminology "*primal*" for any problem where we optimize over $\vartheta \in \Theta$ and "*dual*" for any problem where we optimize over $Q \in \mathbb{P}_{e,f}$.

For the contingent claim $B$, we assume that

$$B \in L^\infty := L^\infty(P).$$

We make this strong assumption because we want results for arbitrary risk aversion parameters $\alpha$. It also has the benefit that our setup fits comfortably into the framework of Delbaen et al. [8]. The measure $P_B$ introduced there via $dP_B := \mathrm{const.} e^{\alpha B} \, dP$ has the same $L^p$ spaces as $P$, and our space $\Theta$ here is the space $\Theta_2$ from [8].

With the above choices, the processes $V^B(x)$ and $C(x, B)$ are well defined for any bounded adapted process $x$ and we get the *exponential utility indifference value process* as

$$(1.2) \quad \begin{aligned} C_t(B) &= \frac{1}{\alpha} \log \frac{V_t^B(0)}{V_t^0(0)} \\ &= \frac{1}{\alpha} \log \left( \operatorname*{ess\,sup}_{\vartheta \in \Theta} E[-e^{-\alpha(G_{t,T}(\vartheta) - B)} | \mathcal{F}_t] \Big/ \operatorname*{ess\,sup}_{\vartheta \in \Theta} E[-e^{-\alpha G_{t,T}(\vartheta)} | \mathcal{F}_t] \right), \end{aligned}$$

independently of the initial capital $x_t$ at time $t$. This yields $C(B)$ in terms of the solutions of two primal problems, but it will be more useful to rewrite this in terms of just one optimization problem. To that end, we introduce the process

$$(1.3) \quad \widetilde{V}_t^B := \operatorname*{ess\,inf}_{Q \in \mathbb{P}_{e,f}} E_Q \left[ \frac{1}{\alpha} \log \frac{Z_T^Q}{Z_t^Q} - B \Big| \mathcal{F}_t \right], \qquad 0 \le t \le T,$$

where $Z^Q$ denotes the density process of $Q$ with respect to $P$. For $B = 0$, $\widetilde{V}^0$ is the dynamic value process associated to the problem of finding the *minimal entropy martingale measure*

$$Q^E := \operatorname*{arg\,min}_{Q \in \mathbb{P}_{e,f}} H(Q|P).$$

In the same way, $\widetilde{V}^B$ is the dynamic value process corresponding to the problem of finding

$$Q^{E,B} := \operatorname*{arg\,min}_{Q \in \mathbb{P}_{e,f}} H(Q|P_B),$$

where $P_B$ is the measure with density $\mathrm{const.} e^{\alpha B}$ with respect to $P$.

PROPOSITION 1. *If $Z^E := Z^{Q^E}$ denotes the density process of $Q^E$ with respect to $P$, then*

$$(1.4) \quad Z_T^E = c_E \exp(G_T(\vartheta^E))$$

*for some constant $c_E \in (0, \infty)$ and some $\vartheta^E \in \Theta$, and*

$$(1.5) \quad \alpha \widetilde{V}_t^0 = E_{Q^E} \left[ \log \frac{Z_T^E}{Z_t^E} \Big| \mathcal{F}_t \right] = \log c_E + G_t(\vartheta^E) - \log Z_t^E, \qquad 0 \le t \le T.$$

*A completely analogous result holds for $Z^{Q^{E,B}}$ and $\widetilde{V}^B$.*

PROOF. The representation (1.4) is well known from [16] and [19]; see Theorem 2.1 of [23] who also prove that the integrand $\vartheta^E$ is in $\Theta$. (1.5) follows from Proposition 4.1 of [23] and (1.4), and the last assertion is obtained by rewriting everything under $P_B$ instead of $P$. □

The next result provides the link between the primal and dual processes $V^B$ and $\widetilde{V}^B$. This is a dynamic version of the results in [8].

PROPOSITION 2. *For fixed $\alpha$, the processes $V^B(0)$ (for $x \equiv 0$) and $\widetilde{V}^B$ are related by*

$$(1.6) \qquad V_t^B(0) = -\exp(-\alpha \widetilde{V}_t^B).$$

*As a consequence, the utility indifference value can be rewritten as*

$$\begin{aligned}
C_t(B) &= -\widetilde{V}_t^B + \widetilde{V}_t^0 \\
&= \operatorname*{ess\,sup}_{Q \in \mathbb{P}_{e,f}} E_Q\left[B - \frac{1}{\alpha}\log\frac{Z_T^Q}{Z_t^Q}\bigg|\mathcal{F}_t\right] - \operatorname*{ess\,sup}_{Q \in \mathbb{P}_{e,f}} E_Q\left[-\frac{1}{\alpha}\log\frac{Z_T^Q}{Z_t^Q}\bigg|\mathcal{F}_t\right].
\end{aligned}$$
(1.7)

PROOF. For $t = 0$, (1.7) is just Theorem 2.2 in [8] whose assumption (2.13) has been shown to be superfluous by Kabanov and Stricker [23]. For general $t \in [0, T]$, the argument is completely analogous; it uses Proposition 1 with $Q^{E,B}$ instead of $Q^E$. □

The representation in (1.7) gives $C(B)$ in terms of the solutions of two dual problems. The desired representation via one single primal problem follows via Proposition 1.

PROPOSITION 3. *The exponential utility indifference value process can be written as*

$$(1.8) \qquad C_t(B) = \frac{1}{\alpha} \log \operatorname*{ess\,inf}_{\vartheta \in \Theta} E_{Q^E}[e^{\alpha(B - G_{t,T}(\vartheta))}|\mathcal{F}_t], \qquad 0 \leq t \leq T.$$

PROOF. If we define the process $\bar{Z} := c_E \exp(G(\vartheta^E))$, then (1.6) and (1.5) tell us that $-V^0(0) = \exp(-\alpha \widetilde{V}^0) = -Z^E/\bar{Z}$, and $\bar{Z}_T = Z_T^E$ by (1.4). Hence, (1.2) yields

$$\begin{aligned}
e^{\alpha C_t(B)} &= -V_t^B(0)\frac{\bar{Z}_t}{Z_t^E} \\
&= -\operatorname*{ess\,sup}_{\vartheta \in \Theta} E\left[-e^{-\alpha(G_{t,T}(\vartheta) - B)}\frac{Z_T^E}{Z_t^E}\frac{\bar{Z}_t}{\bar{Z}_T}\bigg|\mathcal{F}_t\right] \\
&= \operatorname*{ess\,inf}_{\vartheta \in \Theta} E_{Q^E}[e^{\alpha(B - G_{t,T}(\vartheta)) - G_{t,T}(\vartheta^E)}|\mathcal{F}_t].
\end{aligned}$$

Because $\vartheta \mapsto \vartheta' := \vartheta + \vartheta^E/\alpha$ is a bijection from $\Theta$ onto itself, the assertion follows. $\square$

REMARK. Proposition 2 shows that our definition via (1.1) of the utility indifference value *process* agrees with that used in Rouge and El Karoui [33]; see the proof of their Theorem 5.1. Proposition 3 is important because it expresses the utility indifference value process $C(B)$ as the dynamic value process of a standard problem of utility maximization with a random endowment $B$, formulated under the minimal entropy martingale measure $Q^E$. This provides in Section 4 below the link between our dynamic description of $C(B)$ and the recent results of Hu, Imkeller and Müller [21].

**2. Elementary properties of the indifference value.** In this section we list some properties of the exponential utility indifference valuation. These are static properties in the sense that we consider $C_t(B)$ for some fixed $t \in [0,T]$. Our main point is that Propositions 1–3 allow us to extend results known for $t=0$ very easily to arbitrary $t \in [0,T]$. To indicate the dependence on the risk aversion parameter $\alpha$ as well, we write $C_t(B;\alpha)$. Since $P$ is fixed, we write $L^\infty(\mathcal{F}_t)$ for $L^\infty(\mathcal{F}_t, P)$.

PROPOSITION 4. *For fixed $t \in [0,T]$ and $\alpha \in (0,\infty)$, the mapping $B \mapsto C_t(B;\alpha)$ has the following properties:*

(P0) *It maps $L^\infty(\mathcal{F}_T)$ into $L^\infty(\mathcal{F}_t)$, and we have $-\|B\|_\infty \leq C_t(B;\alpha) \leq +\|B\|_\infty$.*
(P1) *It is increasing in $B$: If $B \leq B'$, then $C_t(B;\alpha) \leq C_t(B';\alpha)$.*
(P2) *It is $\mathcal{F}_t$-measurably convex in $B$: we have $C_t(\lambda B + (1-\lambda)B';\alpha) \leq \lambda C_t(B;\alpha) + (1-\lambda)C_t(B';\alpha)$ for any $\lambda \in L^0(\mathcal{F}_t)$ with values in $[0,1]$ and any $B, B' \in L^\infty(\mathcal{F}_T)$.*
(P3) *It is translation-invariant with respect to $L^\infty(\mathcal{F}_t)$ in the sense that we have $C_t(B + x_t;\alpha) = C_t(B;\alpha) + x_t$ for any $x_t \in L^\infty(\mathcal{F}_t)$.*

PROOF. Since $U(z-B) = U(z)e^{\alpha B}$, (P0) is obtained by using the definition of $C_t(B;\alpha)$ via (1.1). (P1)–(P3) follow from the representation (1.7) in Proposition 2 because each functional in the definition (1.3) of $\widetilde{V}_t^B$ has the claimed properties. $\square$

REMARK. In view of Proposition 4, we might call $B \mapsto C_t(B;\alpha)$ a *convex monetary utility functional* from $L^\infty(\mathcal{F}_T)$ to $L^\infty(\mathcal{F}_t)$, because the mapping $B \mapsto C_t(-B;\alpha)$ satisfies the obvious generalizations of the axioms for a convex measure of risk as introduced in [13]; see also [5] for such a suggestion.

While we expect to obtain (P0)–(P2) for $C(B)$ with any reasonable utility function $U$, the next properties are linked to the exponential case.

PROPOSITION 5. *For fixed $t \in [0,T]$, the mapping $B \mapsto C_t(B;\alpha)$ has the following properties:*

(P4) *It does not depend on the initial capital $x_t$ in the definition (1.1).*
(P5) *It is volume-scaling in the sense that $C_t(\beta B; \alpha) = \beta C_t(B; \beta\alpha)$ for any $\beta \in (0, \infty)$.*
(P6) *It is increasing in the risk-aversion $\alpha$: If $\alpha \leq \alpha'$, then $C_t(B;\alpha) \leq C_t(B;\alpha')$.*
(P7) *It satisfies $C_t(\gamma B; \alpha) \leq \gamma C_t(B;\alpha)$ for $\gamma \in [0,1]$ and $C_t(\gamma B;\alpha) \geq \gamma C_t(B;\alpha)$ for $\gamma \in [1,\infty)$.*

PROOF. (P4) is obvious, (P7) follows directly from (P5) and (P6), and these are proved via the representation (1.8) in Proposition 3; (P5) uses that $\Theta$ is a cone, (P6) uses Jensen's inequality. □

The preceding results are in no way original; they go back to Rouge and El Karoui [33] and Becherer [3] who formulated and proved them for $t = 0$. These authors also gave asymptotic results for large and small risk aversions ($\alpha \nearrow \infty$ and $\alpha \searrow 0$) and we shall prove below versions of those results for arbitrary $t \in [0,T]$ with the help of a description of the process $(C_t(B;\alpha))_{0 \leq t \leq T}$ via a backward stochastic differential equation. Before we embark on that aspect, however, we give two more properties of $C_t(B)$. The first says that anything which is attainable at zero cost by self-financing trading between $t$ and $T$ has zero value and does not affect the valuation of $B$; the second says that $C_t(B)$ always lies in the interval of arbitrage-free prices for $B$. Such results for $t = 0$ have already been given by Rouge and El Karoui [33] and Becherer [3], among others; see also [17].

LEMMA 6. *For any $t \in [0,T]$ and $\alpha \in (0,\infty)$, we have the following:*

(1) *For any $\vartheta \in \Theta$, $C_t(G_{t,T}(\vartheta); \alpha) = 0$ and $C_t(B + G_{t,T}(\vartheta); \alpha) = C_t(B; \alpha)$.*
(2) $\operatorname{ess\,inf}_{Q \in \mathbb{P}_{e,f}} E_Q[B|\mathcal{F}_t] \leq C_t(B;\alpha) \leq \operatorname{ess\,sup}_{Q \in \mathbb{P}_{e,f}} E_Q[B|\mathcal{F}_t]$.

PROOF. (1) Since $G(\vartheta)$ is a $Q$-martingale for any $Q \in \mathbb{P}_{e,f}$, this is immediate from (1.7).
(2) We know from (1.7) and (1.3) that

$$C_t(B;\alpha) = \operatorname*{ess\,sup}_{Q \in \mathbb{P}_{e,f}} \left( E_Q[B|\mathcal{F}_t] - \frac{1}{\alpha} \left( E_Q\left[ \log \frac{Z_T^Q}{Z_t^Q} \Big| \mathcal{F}_t \right] - \alpha \widetilde{V}_t^0 \right) \right).$$

By the definition of $\widetilde{V}^0$ in (1.3), the term in the inner brackets is always nonnegative, and it equals zero for $Q = Q^E$ by Proposition 1. The first fact gives the upper bound in (2), the second one the lower bound. □

**3. A comparison theorem and some results for a BSDE.** This section studies a family of backward stochastic differential equations (BSDEs) that play an important role in a dynamic description of the exponential utility indifference value. We work on a filtered probability space $(\Omega, \mathcal{F}, \mathbb{F}, R)$ and we assume throughout this section that the filtration $\mathbb{F}$ is *continuous*, that is, all local martingales are continuous.

We fix a (continuous) $\mathbb{R}^d$-valued local $R$-martingale $M$ null at 0 and denote by $BMO[M]$ the space of all $\mathbb{R}^d$-valued predictable $M$-integrable processes $h$ such that $h \cdot M := \int h \, dM$ is in $BMO(R)$, the usual martingale space $BMO$ for the measure $R$. Note that $\langle M \rangle$ is a $(d \times d)$ matrix-valued process.

Let us consider the semimartingale backward equation

$$(3.1) \quad Y_t = Y_0 + \int_0^t \underline{1}^{\text{tr}} d\langle M \rangle_s f(s, Z_s) + \int_0^t g_s \, d\langle L \rangle_s + \int_0^t Z_s \, dM_s + L_t$$

with the boundary condition

$$(3.2) \quad Y_T = B,$$

where $\underline{1} := (1 \ldots 1)^{\text{tr}} \in \mathbb{R}^d$, $f : \Omega \times [0,T] \times \mathbb{R}^d \to \mathbb{R}^d$ is $\mathcal{P} \times \mathcal{B}(\mathbb{R}^d)$-measurable, $g$ is a real-valued predictable process and $B \in L^\infty(\mathcal{F}_T, R)$. We call $(f, g, B)$ the generator of (3.1) and (3.2). A *solution* of (3.1) and (3.2) is a triple $(Y, Z, L)$, where $Y$ is a real-valued special $R$-semimartingale, $Z$ is an $\mathbb{R}^d$-valued predictable $M$-integrable process and $L$ is a real-valued local $R$-martingale strongly $R$-orthogonal to $M$. Sometimes we call $Y$ alone the solution of (3.1) and (3.2), keeping in mind that $Z \cdot M + L$ is the martingale part of $Y$.

Our first result and its subsequent applications show the importance of $BMO$-martingales when dealing with BSDEs with quadratic generators; see also [21, 28] or [26].

PROPOSITION 7. *Suppose there are constants $C_f, C_g$ and a predictable process $K \in BMO[M]$ such that*

$$C_f \int Z_s^{\text{tr}} d\langle M \rangle_s Z_s + \int K_s^{\text{tr}} d\langle M \rangle_s K_s - \int |\underline{1}^{\text{tr}} d\langle M \rangle_s f(s, Z_s)|$$
*is an increasing process for any $\mathbb{R}^d$-valued predictable $M$-integrable $Z$,*
(3.3)

$$(3.4) \quad |g_t| \leq C_g, \qquad R\text{-a.s., for each } t \in [0,T].$$

*Then the martingale part of any bounded solution of (3.1) and (3.2) is in $BMO(R)$.*

PROOF. Let $Y$ be a solution of (3.1) and (3.2) and $c > 0$ a constant such that

(3.5) $$|Y_t| \leq c, \quad R\text{-a.s., for each } t \in [0,T].$$

Applying Itô's formula between a stopping time $\tau$ and $T$ and using (3.5) yields

(3.6) $$\begin{aligned} e^{|\beta|c} &\geq e^{\beta Y_T} - e^{\beta Y_\tau} \\ &= \frac{\beta^2}{2} \int_\tau^T e^{\beta Y_s} Z_s^{\mathrm{tr}} d\langle M \rangle_s Z_s + \frac{\beta^2}{2} \int_\tau^T e^{\beta Y_s} d\langle L \rangle_s \\ &\quad + \beta \int_\tau^T e^{\beta Y_s} \underline{1}^{\mathrm{tr}} d\langle M \rangle_s f(s, Z_s) + \beta \int_\tau^T e^{\beta Y_s} g_s d\langle L \rangle_s \\ &\quad + \beta \int_\tau^T e^{\beta Y_s} Z_s dM_s + \beta \int_\tau^T e^{\beta Y_s} dL_s, \end{aligned}$$

where $\beta \in \mathbb{R}$ is a constant yet to be determined.

If $Z \cdot M$ and $L$ are true $R$-martingales, taking conditional expectations in (3.6) gives

$$\frac{\beta^2}{2} E_R\left[\int_\tau^T e^{\beta Y_s} Z_s^{\mathrm{tr}} d\langle M \rangle_s Z_s \Big| \mathcal{F}_\tau\right] + \frac{\beta^2}{2} E_R\left[\int_\tau^T e^{\beta Y_s} d\langle L \rangle_s \Big| \mathcal{F}_\tau\right]$$
$$\leq e^{|\beta|c} + |\beta| E_R\left[\int_\tau^T e^{\beta Y_s} |\underline{1}^{\mathrm{tr}} d\langle M \rangle_s f(s, Z_s)| \Big| \mathcal{F}_\tau\right]$$
$$+ |\beta| E_R\left[\int_\tau^T e^{\beta Y_s} |g_s| d\langle L \rangle_s \Big| \mathcal{F}_\tau\right].$$

Using the conditions (3.3) and (3.4), we can rewrite this estimate as

(3.7) $$\begin{aligned} &\left(\frac{\beta^2}{2} - |\beta| C_f\right) E_R\left[\int_\tau^T e^{\beta Y_s} Z_s^{\mathrm{tr}} d\langle M \rangle_s Z_s \Big| \mathcal{F}_\tau\right] \\ &\quad + \left(\frac{\beta^2}{2} - |\beta| C_g\right) E_R\left[\int_\tau^T e^{\beta Y_s} d\langle L \rangle_s \Big| \mathcal{F}_\tau\right] \\ &\leq e^{|\beta|c} + |\beta| E_R\left[\int_\tau^T e^{\beta Y_s} K_s^{\mathrm{tr}} d\langle M \rangle_s K_s \Big| \mathcal{F}_\tau\right] \\ &\leq e^{|\beta|c}(1 + |\beta| \|K \cdot M\|_{\mathrm{BMO}(R)}^2). \end{aligned}$$

For $\beta := 4\overline{C} := 4\max(C_f, C_g) > 0$, we obtain from (3.7) that

$$4\overline{C}^2 \left( E_R\left[\int_\tau^T e^{\beta Y_s} Z_s^{\mathrm{tr}} d\langle M \rangle_s Z_s \Big| \mathcal{F}_\tau\right] + E_R\left[\int_\tau^T e^{\beta Y_s} d\langle L \rangle_s \Big| \mathcal{F}_\tau\right] \right)$$
$$\leq e^{4\overline{C}c}(1 + 4\overline{C}\|K \cdot M\|_{\mathrm{BMO}(R)}^2),$$

and if we use (3.5) to write $e^{\beta Y_s} \geq e^{-|\beta|c} = e^{-4\overline{C}c}$, we finally get

(3.8)
$$E\left[\int_\tau^T Z_s^{\text{tr}} d\langle M\rangle_s Z_s \Big| \mathcal{F}_\tau\right] + E[\langle L\rangle_T - \langle L\rangle_\tau | \mathcal{F}_\tau]$$
$$\leq \frac{e^{8\overline{C}c}(1 + 4\overline{C}\|K \cdot M\|^2_{BMO(R)})}{4\overline{C}^2},$$

$R$-a.s. for any stopping time $\tau$. Hence, $Z \cdot M$ and $L$ are in $BMO(R)$.

For general $Z \cdot M$ and $L$, we stop at $\tau_n$ and apply the above argument with $T$ replaced by $\tau_n$ to get (3.8) also with $T$ replaced by $\tau_n$. Letting $n \to \infty$ then completes the proof. □

We are now in a position to give a comparison theorem for the BSDE (3.1). Although we need this result only for $f \equiv 0$, we formulate and prove it in general.

THEOREM 8. *Suppose the generators $(f^i, g^i, B^i)$, $i = 1, 2$, satisfy the assumptions of Proposition 7, and $Y^i$, $i = 1, 2$, are corresponding bounded solutions of (3.1) and (3.2). (In particular, we assume here the existence of these solutions.) Suppose also that $B^1 \geq B^2$, $R$-a.s.; that the process $\int \mathbf{1}^{\text{tr}} d\langle M\rangle_s (f^1(s, Z_s) - f^2(s, Z_s))$ is decreasing for any $Z \in BMO[M]$; that $g^1 \leq g^2$ $R \otimes \langle L\rangle$-a.e.; and that either $f^1$ or $f^2$ satisfies the following condition:*

*For any $Z^1, Z^2 \in BMO[M]$, there exists some $\kappa \in BMO[M]$ such that*
$$\int \mathbf{1}^{\text{tr}} d\langle M\rangle_s (f(s, Z_s^1) - f(s, Z_s^2)) = \int \kappa_s^{\text{tr}} d\langle M\rangle_s (Z_s^1 - Z_s^2).$$
(3.9)
*Then $Y_t^1 \geq Y_t^2$ $R$-a.s. for all $t \in [0, T]$.*

PROOF. By taking differences, we obtain
$$Y_t^1 - Y_t^2 - (Y_0^1 - Y_0^2)$$
$$= \int_0^t \mathbf{1}^{\text{tr}} d\langle M\rangle_s (f^1(s, Z_s^2) - f^2(s, Z_s^2)) + \int_0^t (g_s^1 - g_s^2) d\langle L^1\rangle_s$$
$$+ \int_0^t \mathbf{1}^{\text{tr}} d\langle M\rangle_s (f^1(s, Z_s^1) - f^1(s, Z_s^2)) + \int_0^t g_s^2 d(\langle L^1\rangle_s - \langle L^2\rangle_s)$$
$$+ \int_0^t (Z_s^1 - Z_s^2) dM_s + L_t^1 - L_t^2.$$

Suppose $f^1$ satisfies (3.9). According to Proposition 7, $Z^1 \cdot M$, $Z^2 \cdot M$, $L^1$, $L^2$ are all in $BMO(R)$. Hence, (3.9) and (3.4) imply that
$$N := -\int \kappa_s dM_s - \int g_s^2 d(L_s^1 + L_s^2)$$

is in $BMO(R)$, and so $Q$ defined by $dQ = \mathcal{E}(N)_T \, dR$ is a probability measure equivalent to $R$; see Theorem 2.3 of [24]. If

$$\bar{N} := (Z^1 - Z^2) \cdot M + L^1 - L^2$$

denotes the $R$-martingale part of $Y^1 - Y^2$, (3.9) yields that

$$
\begin{aligned}
Y^1 - Y^2 &- (Y_0^1 - Y_0^2) - \int \underline{1}^{\mathrm{tr}} \, d\langle M \rangle_s (f^1(s, Z_s^2) - f^2(s, Z_s^2)) \\
&\quad - \int (g_s^1 - g_s^2) \, d\langle L^1 \rangle_s \\
&= \int \underline{1}^{\mathrm{tr}} \, d\langle M \rangle_s (f^1(s, Z_s^1) - f^1(s, Z_s^2)) + \int g_s^2 \, d(\langle L^1 \rangle_s - \langle L^2 \rangle_s) + \bar{N} \\
&= \bar{N} + \int \kappa_s^{\mathrm{tr}} \, d\langle M \rangle_s (Z_s^1 - Z_s^2) + \int g_s^2 \, d(\langle L^1 \rangle_s - \langle L^2 \rangle_s) \\
&= \bar{N} - \langle N, \bar{N} \rangle
\end{aligned}
\tag{3.10}
$$

is a local $Q$-martingale by Girsanov's theorem and even in $BMO(Q)$ by Theorem 3.6 of [24], since $\bar{N}$ is in $BMO(R)$ by Proposition 7. Thus, we can use the $Q$-martingale property and the boundary conditions $Y_T^i = B^i$ to obtain from (3.10) that

$$
\begin{aligned}
Y_t^1 - Y_t^2 &\\
&= E_Q\bigg[ B^1 - B^2 - \int_t^T \underline{1}^{\mathrm{tr}} \, d\langle M \rangle_s (f^1(s, Z_s^2) - f^2(s, Z_s^2)) \\
&\qquad\qquad - \int_t^T (g_s^1 - g_s^2) \, d\langle L^1 \rangle_s \bigg| \mathcal{F}_t \bigg],
\end{aligned}
\tag{3.11}
$$

which implies the assertion. $\square$

REMARKS. (1) The assumption (3.3) is a quadratic condition (in $Z$). This becomes more apparent if we use the strong order on increasing processes (where $A \preceq A'$ means that $A' - A$ is increasing) to rewrite (3.3) more compactly as

$$\int |\underline{1}^{\mathrm{tr}} \, d\langle M \rangle_s f(s, Z_s)| \preceq C_f \langle Z \cdot M \rangle + \langle K \cdot M \rangle.$$

(2) For $d = 1$, the BSDE (3.1) and the above conditions on $f$ take a simpler and more familiar form since $\langle M \rangle$ is then a scalar process. The term $\int \underline{1}^{\mathrm{tr}} \, d\langle M \rangle_s \, f(s, Z_s)$ in (3.1) reduces to $\int f(s, Z_s) \, d\langle M \rangle_s$; the condition on $f^1 - f^2$ in Theorem 8 follows if $f^1(t, z) \leq f^2(t, z)$; (3.3) boils down to the quadratic growth condition $|f(t, z)| \leq K_t^2 + C_f z^2$; and (3.9) essentially means that (with $0/0 := 0$)

$$\frac{f(\cdot, Z_\cdot^1) - f(\cdot, Z_\cdot^2)}{Z_\cdot^1 - Z_\cdot^2} \in BMO[M] \qquad \text{for any } Z^1, Z^2 \in BMO[M]. \tag{3.12}$$

Note that this is fulfilled for functionals of the form $f(\omega, t, z) = D_t^0(\omega) + D_t^1(\omega)z + D_t^2(\omega)z^2$ with processes $D^0, D^1$ in $BMO[M]$ and a bounded predictable process $D^2 \geq 0$. Alternatively, (3.12) holds if $f(t, z)$ satisfies a global Lipschitz condition in $z$ and $M$ is in $BMO(R)$.

For later use, we consider the special case of the generator $(0, -\frac{\alpha}{2}, B)$ with $\alpha \in (0, \infty)$ and $B \in L^\infty(R)$. The BSDE (3.1) then takes the form (with $\psi$ replacing $Z$)

$$(3.13) \qquad Y_t = Y_0 - \frac{\alpha}{2} \langle L \rangle_t + \int_0^t \psi_s \, dM_s + L_t,$$

and its solution with final condition $Y_T = B$ is denoted by $(Y^\alpha, \psi^\alpha, L^\alpha)$. We now derive estimates on these quantities as $\alpha$ varies.

LEMMA 9. *For the solutions $(Y^\alpha, \psi^\alpha, L^\alpha)$ of (3.13) and (3.2) with generator $(0, -\frac{\alpha}{2}, B)$, we have*

$$(3.14) \qquad \sup_{\alpha \in (0,\infty)} \|\psi^\alpha \cdot M\|_{BMO(R)} < \infty,$$

$$(3.15) \qquad \sup_{\alpha \in (0,\infty)} \alpha \|L^\alpha\|_{BMO(R)}^2 < \infty.$$

*In particular, this yields*

$$(3.16) \qquad \sup_{\alpha \in (0,\infty)} \|L^\alpha\|_{BMO(R)} < \infty,$$

$$(3.17) \qquad \lim_{\alpha \to \infty} \|L^\alpha\|_{BMO(R)} = 0.$$

PROOF. We go back to the proof of Proposition 7 and note that $C_f = 0$, $K \equiv 0$ in (3.3) and $C_g = \frac{\alpha}{2}$ in (3.4). Hence, we obtain from (3.6) as for (3.7) with $\beta = -1$ and $c = \|B\|_\infty$ that

$$e^{\|B\|_\infty} \geq \frac{1}{2} e^{-\|B\|_\infty} E_R \left[ \int_\tau^T (\psi_u^\alpha)^{\mathrm{tr}} \, d\langle M \rangle_u \, \psi_u^\alpha \Big| \mathcal{F}_\tau \right]$$
$$+ \frac{1+\alpha}{2} e^{-\|B\|_\infty} E_R[\langle L^\alpha \rangle_T - \langle L^\alpha \rangle_\tau | \mathcal{F}_\tau],$$

where we have used in (3.6) $\beta g_s \equiv \frac{\alpha}{2}$ instead of the cruder estimate $\beta g_s \geq -|\beta| C_g = -\frac{\alpha}{2}$. The above estimate yields

$$(3.18) \qquad \|\psi^\alpha \cdot M\|_{BMO(R)}^2 + (1+\alpha) \|L^\alpha\|_{BMO(R)}^2 \leq 2 e^{2\|B\|_\infty}$$

uniformly for all $\alpha \in (0, \infty)$.

Thus, we obtain (3.14) and (3.15), and (3.16) and (3.17) then follow immediately. □

REMARK. One can also deduce (3.15)–(3.17) by taking conditional expectations directly in (3.13). We have chosen the above argument since it gives (3.14) at the same time.

PROPOSITION 10. *The solution $Y^\alpha$ of (3.13) and (3.2) is locally Lipschitz-continuous with respect to $\alpha$, uniformly in $t$: For any $\gamma > 0$, there is a constant $K_\gamma$ depending only on $\gamma$ such that*

$$(3.19) \quad \sup_{0 \leq t \leq T} |Y_t^\alpha - Y_t^{\alpha'}| \leq K_\gamma |\alpha - \alpha'| \quad \text{for all } \alpha, \alpha' \in (0, \gamma].$$

PROOF. We go back to the proof of Theorem 8 with the two generators $(0, -\frac{\alpha}{2}, B)$ and $(0, -\frac{\alpha'}{2}, B)$. Then (3.11) yields

$$(3.20) \quad Y_t^\alpha - Y_t^{\alpha'} = \frac{\alpha - \alpha'}{2} E_Q[\langle L^\alpha \rangle_T - \langle L^\alpha \rangle_t | \mathcal{F}_t],$$

where $Q$ is now given by

$$dQ = \mathcal{E}\left(\frac{\alpha'}{2}(L^\alpha + L^{\alpha'})\right)_T dR =: \mathcal{E}\left(\frac{\alpha'}{2} L(\alpha, \alpha')\right)_T dR =: Z_T(\alpha, \alpha') \, dR.$$

Due to (3.16), we have

$$(3.21) \quad \begin{aligned} &\sup_{\alpha, \alpha' \in (0, \gamma]} \left\| \frac{\alpha'}{2} L(\alpha, \alpha') \right\|_{BMO(R)} \\ &\leq \frac{\gamma}{2} \sup_{\alpha, \alpha' \in (0, \infty)} (\|L^\alpha\|_{BMO(R)} + \|L^{\alpha'}\|_{BMO(R)}) < \infty. \end{aligned}$$

By Theorem 3.1 of [24], $Z(\alpha, \alpha')$ therefore satisfies the reverse Hölder inequality $\mathcal{R}_p(R)$ for some $p \in (1, \infty)$, that is,

$$\sup_{0 \leq t \leq T} E_R\left[\left(\frac{Z_T(\alpha, \alpha')}{Z_t(\alpha, \alpha')}\right)^p \bigg| \mathcal{F}_t \right] \leq (c_p)^p$$

for a constant $c_p$; this holds uniformly for all $\alpha, \alpha' \in (0, \gamma]$ since (3.21) is also uniform in those $\alpha, \alpha'$. Moreover, the energy inequalities (see [24], page 28) yield

$$(3.22) \quad \sup_{0 \leq t \leq T} E_R[(\langle L^\alpha \rangle_T - \langle L^\alpha \rangle_t)^n | \mathcal{F}_t] \leq n! \|L^\alpha\|_{BMO(R)}^{2n} \quad \text{for all } n \in \mathbb{N}.$$

So if we choose $n$ with $\frac{n}{n-1} \leq p$, Bayes' rule and Hölder's inequality give

$$\sup_{0 \leq t \leq T} E_Q[\langle L^\alpha \rangle_T - \langle L^\alpha \rangle_t | \mathcal{F}_t] = \sup_{0 \leq t \leq T} E_R\left[\frac{Z_T(\alpha, \alpha')}{Z_t(\alpha, \alpha')}(\langle L^\alpha \rangle_T - \langle L^\alpha \rangle_t) \bigg| \mathcal{F}_t \right]$$

$$\leq c_{n/(n-1)} \sup_{0 \leq t \leq T} (E_R[(\langle L^\alpha \rangle_T - \langle L^\alpha \rangle_t)^n | \mathcal{F}_t])^{1/n}.$$

Combining this with (3.20), (3.22) and (3.21) yields (3.19). □

A closer look at the proof of Theorem 8 shows that we can also write down a quasi-explicit expression for $Y^\alpha$.

PROPOSITION 11. *The solution* $(Y^\alpha, \psi^\alpha, L^\alpha)$ *of* (3.13) *and* (3.2) *with generator* $(0, -\frac{\alpha}{2}, B)$ *can be represented as follows: If we define the measure* $Q^\alpha$ *by* $dQ^\alpha := \tilde{Z}_T^\alpha dR := \mathcal{E}(\frac{\alpha}{2} L^\alpha)_T dR$, *then*

$$
\begin{aligned}
Y_t^\alpha &= E_{Q^\alpha}[B|\mathcal{F}_t] \\
&= E_R\bigg[\frac{\mathcal{E}((\alpha/2)L^\alpha)_T}{\mathcal{E}((\alpha/2)L^\alpha)_t} B \bigg| \mathcal{F}_t\bigg], \qquad R\text{-}a.s. \text{ for each } t \in [0,T],
\end{aligned}
\tag{3.23}
$$

*and* $\psi^\alpha$ *is a predictable density of* $\langle Y^\alpha, M\rangle$ *with respect to* $\langle M\rangle$, *that is,* $d\langle Y^\alpha, M\rangle = d\langle M\rangle \psi^\alpha$.

PROOF. For the two generators $(0, -\frac{\alpha}{2}, B)$ and $(0, -\frac{\alpha}{2}, 0)$ with corresponding solutions $(Y^\alpha, \psi^\alpha, L^\alpha)$ and $(0, 0, 0)$, the martingale $N$ in the proof of Theorem 8 reduces to $\frac{\alpha}{2} L^\alpha$. Hence, (3.23) follows from (3.11), and the second assertion then from the BSDE (3.13). □

Note that the representation (3.23) of $Y^\alpha$ is not as simple as it may appear, because the measure $Q^\alpha$ still involves the component $L^\alpha$ from the solution triple $(Y^\alpha, \psi^\alpha, L^\alpha)$. Since this depends on $B$ via the final condition (3.2), (3.23) is, in particular, not linear in $B$ in general.

**4. Dynamic description of the utility indifference value.** In this section we study the dynamic behavior of the exponential utility indifference value over time. We characterize the process $C(B; \alpha)$ as the unique solution of a BSDE in a general continuous filtration which need not be generated by a Brownian motion, thus extending earlier results by Rouge and El Karoui [33]. Given the characterization of $C(B; \alpha)$ in Proposition 3, we can also view our BSDE as a generalization of the one obtained independently by Hu, Imkeller and Müller [21]. Finally, our BSDE is also a continuous-time analogue of the recursive description in Theorem 5 of [30], obtained in a particular discrete-time setting.

To prove existence and uniqueness of a solution to their BSDEs, Rouge and El Karoui [33] and Hu, Imkeller and Müller [21] used results of Kobylanski [25] on existence and comparison for quadratic BSDEs driven by a Brownian motion. But for BSDEs with quadratic generators and driven by martingales, there are no general results similar to those of Kobylanski [25]. Chitashvili [7] and El Karoui and Huang [11] established the well-posedness of BSDEs driven by martingales if the generators satisfy global Lipschitz conditions, but this is too restrictive for our needs. We prove here existence of a solution by directly showing that $C(B;\alpha)$ satisfies a quadratic BSDE, and we use the comparison theorem from Section 3 to obtain uniqueness.

We start by recalling from Proposition 3 that the exponential utility indifference value process $C(B;\alpha)$ can be represented as

$$(1.8) \quad C_t(B;\alpha) = \frac{1}{\alpha} \log \operatorname*{ess\,inf}_{\vartheta \in \Theta} E_{Q^E}[e^{\alpha(B - G_{t,T}(\vartheta))} | \mathcal{F}_t], \qquad 0 \leq t \leq T.$$

This shows that $e^{\alpha C(B;\alpha)}$ is the dynamic value process of the stochastic control problem

$$(4.1) \qquad \text{minimize } E_{Q^E}[e^{\alpha(B - G_T(\vartheta))}] \qquad \text{over all } \vartheta \in \Theta.$$

Using similar arguments as in [8], one can show that an optimal strategy $\vartheta^* \in \Theta$ for (4.1) exists. The *martingale optimality principle* takes here the following form.

PROPOSITION 12. *Suppose that $S$ is locally bounded, $\mathbb{P}_{e,f} \neq \varnothing$ and $B \in L^\infty$. Fix $\alpha > 0$.*

(1) *There exists an RCLL process $J^B = (J^B_t)_{0 \leq t \leq T}$ such that, for each $t \in [0,T]$,*

$$(4.2) \qquad J^B_t = \operatorname*{ess\,inf}_{\vartheta \in \Theta} E_{Q^E}[e^{\alpha(B - G_{t,T}(\vartheta))} | \mathcal{F}_t], \qquad P\text{-a.s.}$$

*$J^B$ is the largest RCLL process $J$ with $J_T = e^{\alpha B}$, $P$-a.s. such that $Je^{-\alpha G(\vartheta)}$ is a $Q^E$-submartingale for each $\vartheta \in \Theta$.*

(2) *The following properties are equivalent:*

   (a) *$\vartheta^* \in \Theta$ is optimal for (4.1), that is, $J^B_0 = E_{Q^E}[e^{\alpha(B - G_T(\vartheta^*))}]$.*
   (b) *$\vartheta^* \in \Theta$ is optimal for all conditional criteria, that is,*

   $$J^B_t = E_{Q^E}[e^{\alpha(B - G_{t,T}(\vartheta^*))} | \mathcal{F}_t], \qquad P\text{-a.s., for each } t \in [0,T].$$

   (c) *The process $J^B e^{-\alpha G(\vartheta^*)}$ with $\vartheta^* \in \Theta$ is a $Q^E$-martingale.*

(3) *Due to (4.2) and (1.8), we can and do choose $\frac{1}{\alpha} \log J^B$ as an RCLL version for $C(B;\alpha)$. For any stopping times $\sigma \leq \tau \leq T$, we then have the dynamic programming equation*

$$(4.3) \qquad C_\sigma(B;\alpha) = \frac{1}{\alpha} \log \operatorname*{ess\,inf}_{\vartheta \in \Theta} E_{Q^E}[e^{\alpha(C_\tau(B;\alpha) - G_{\sigma,\tau}(\vartheta))} | \mathcal{F}_\sigma], \qquad P\text{-a.s.}$$

PROOF. This is a standard argument like in [12] or [27] and therefore omitted.

□

Because we have an optimal strategy $\vartheta^* \in \Theta$, Proposition 12 yields that

$$C(B;\alpha) = \frac{1}{\alpha}\log J^B = \frac{1}{\alpha}\log(J^B e^{-\alpha G(\vartheta^*)}) + G(\vartheta^*)$$

is a $Q^E$-supermartingale; see Proposition 6 of [30] for an analogous result in a particular discrete-time setting. To obtain more structure for $C(B;\alpha)$, we now assume that

$$\mathbb{F} \text{ is continuous;}$$

this implies, in particular, that $S$ is continuous. The Doob–Meyer decomposition of $C(B;\alpha)$ is

$$C(B;\alpha) = C_0(B;\alpha) + M^B(\alpha) - A^B(\alpha) \qquad \text{under } Q^E,$$

where $M^B(\alpha) \in \mathcal{M}_{0,\mathrm{loc}}(Q^E)$ and $A^B(\alpha)$ is adapted, continuous and increasing. Using the Galtchouk–Kunita–Watanabe decomposition for $M^B(\alpha)$ with respect to $S$ under $Q^E$, we get

(4.4) $$C(B;\alpha) = C_0(B;\alpha) - A^B(\alpha) + \int \varphi^B(\alpha)\,dS + m^B(\alpha)$$

with $m^B(\alpha) \in \mathcal{M}_{0,\mathrm{loc}}(Q^E)$ satisfying $\langle m^B(\alpha), S\rangle = 0$.

THEOREM 13. *Suppose that $\mathbb{P}_{e,f} \neq \varnothing$, $B \in L^\infty$ and $\mathbb{F}$ is continuous. Then the exponential utility indifference value process $C(B;\alpha)$ is the unique bounded solution of the following semimartingale backward equation under the minimal entropy martingale measure $Q^E$:*

(4.5) $$Y_t = Y_0 - \frac{\alpha}{2}\langle L\rangle_t + \int_0^t \psi_s\,dS_s + L_t$$

*with the boundary condition*

(3.2) $$Y_T = B.$$

*["Under $Q^E$" means that, in the solution triple $(Y,\psi,L)$, the process $L$ is a local $Q^E$-martingale strongly $Q^E$-orthogonal to $S$.] Moreover, $\psi \cdot S$ and $L$ are both in $BMO(Q^E)$.*

PROOF. (1) We first show that $C(B;\alpha)$ satisfies (4.5) and (3.2). Applying Itô's formula for $Z^{(\vartheta)} := e^{\alpha(C(B;\alpha)-G(\vartheta))}$ and omitting the index $\alpha$, we have from (4.4)

$$Z_t^{(\vartheta)} = Z_0^{(\vartheta)}$$
$$+ \alpha \int_0^t Z_s^{(\vartheta)}\,d\left(-A_s^B + \frac{\alpha}{2}\int_0^s (\varphi_u^B - \vartheta_u)^{\mathrm{tr}}\,d\langle S\rangle_u\,(\varphi_u^B - \vartheta_u) + \frac{\alpha}{2}\langle m^B\rangle_s\right)$$
$$+ \text{local } Q^E\text{-martingale.}$$
(4.6)

By parts (1) and (2) of Proposition 12, $Z^{(\vartheta)}$ is a $Q^E$-submartingale for any $\vartheta \in \Theta$ and a $Q^E$-martingale for the optimal strategy $\vartheta^*$. Since $Z^{(\vartheta)} > 0$, this implies by (4.6) that

$$-A^B + \frac{\alpha}{2}\int(\varphi^B - \vartheta)^{\mathrm{tr}}\,d\langle S\rangle\,(\varphi^B - \vartheta) + \frac{\alpha}{2}\langle m^B\rangle \quad \text{is increasing}$$

for any $\vartheta \in \Theta$ and vanishes for $\vartheta^*$. Hence, it follows that

$$\begin{aligned}
A^B &= \operatorname*{ess\,inf}_{\vartheta\in\Theta}\left(\frac{\alpha}{2}\int(\varphi^B - \vartheta)^{\mathrm{tr}}\,d\langle S\rangle\,(\varphi^B - \vartheta) + \frac{\alpha}{2}\langle m^B\rangle\right) \\
&= \frac{\alpha}{2}\langle m^B\rangle + \frac{\alpha}{2}\operatorname*{ess\,inf}_{\vartheta\in\Theta}\int(\varphi^B - \vartheta)^{\mathrm{tr}}\,d\langle S\rangle\,(\varphi^B - \vartheta),
\end{aligned} \tag{4.7}$$

where we can take the ess inf with respect to the strong order. To prove that

$$A^B = \frac{\alpha}{2}\langle m^B\rangle, \tag{4.8}$$

we define the stopping times $\tau_n := \inf\{t \geq 0 | |G_t(\varphi^B)| \geq n\}$. Then $\tau_n \nearrow T$ stationarily, $P$-a.s., and $\vartheta^n := \varphi^B I_{]\!]0,\tau_n]\!]}$ is in $\Theta$ for every $n$. Hence, we get, for any $t \leq T$, that

$$\operatorname*{ess\,inf}_{\vartheta\in\Theta}\frac{\alpha}{2}\int_0^t(\varphi_s^B - \vartheta_s)^{\mathrm{tr}}\,d\langle S\rangle_s\,(\varphi_s^B - \vartheta_s) \leq \frac{\alpha}{2}\int_0^t(\varphi_s^B - \vartheta_s^n)^{\mathrm{tr}}\,d\langle S\rangle_s\,(\varphi_s^B - \vartheta_s^n)$$

$$= \frac{\alpha}{2}\int_{\tau_n}^{t\vee\tau_n}(\varphi_s^B)^{\mathrm{tr}}\,d\langle S\rangle_s\,\varphi_s^B \longrightarrow 0$$

as $n \to \infty$, which implies (4.8). Combining this with (4.4) shows that $C(B;\alpha)$ indeed satisfies (4.5), and it is clear that we also have the boundary condition $C_T(B;\alpha) = B$. The *BMO* property of $\psi \cdot S$ and $L$ follows from Proposition 7, applied with the pair $(M,R) = (S,Q^E)$.

(2) We already know from Proposition 4 that $C(B;\alpha)$ is bounded by $\|B\|_\infty$. The uniqueness of a bounded solution of (4.5) and (3.2) follows from the comparison in Theorem 8, applied with the pair $(M,R) = (S,Q^E)$. □

REMARKS. (1) In comparison to the work of Rouge and El Karoui [33] and Hu, Imkeller and Müller [21], our BSDE result in Theorem 13 is at the same time more and less general. We are able to work in a general continuous filtration, but we have so far not included any constraints in our strategies. For the case where $dS_t = \sigma_t\,dW_t^*$ under an equivalent martingale measure $Q^*$ and $\mathbb{F}$ is generated by a Brownian motion, our BSDE (4.5) can be rewritten as

$$dY_t = -\frac{\alpha}{2}|\Pi_t z_t|^2\,dt + z_t\,dW_t^E \qquad \text{under } Q^E,$$

where $\Pi_t$ denotes the projection on $\ker(\sigma_t) = (\operatorname{range}(\sigma_t^{\mathrm{tr}}))^\perp$. This agrees with the BSDEs of Rouge and El Karoui [33] and Hu, Imkeller and Müller [21] in that particular case.

(2) One advantage of our approach is that even in a Brownian filtration, we need not invoke general results on quadratic BSDEs. This allows us to avoid restrictive assumptions (like boundedness) on the coefficients of our model. In fact, our only requirement is the natural condition that the minimal entropy martingale measure $Q^E$ exists.

(3) The proof of Theorem 13 shows, in particular, that the value of the infimum in (4.7) is obtained by choosing $\vartheta = \varphi^B$. Because we already know that an optimal strategy $\vartheta^* \in \Theta$ exists, we conclude that $\vartheta^* = \varphi^B$, and, in particular, that $\varphi^B$ is in $\Theta$. Moreover, we also see from (4.4) that the $\psi$-component of the solution to the BSDE (4.5) is given by the optimal strategy $\vartheta^*$ for the utility maximization problem (4.1).

(4) If we only assume that $S$ is continuous while the filtration is general, we can still show that $C(B;\alpha)$ satisfies the semimartingale backward equation

$$(4.9) \quad Y_t = Y_0 - \frac{1}{\alpha}\left(\sum_{0<s\leq t}(e^{\alpha \Delta Y_s} - \alpha \Delta Y_s - 1)\right)^p - \frac{\alpha}{2}\langle L\rangle_t + \int_0^t \psi_s\, dS_s + L_t$$

with boundary condition $Y_T = B$, where $A^p$ denotes the dual predictable projection of a locally integrable increasing process $A$. We do not have a comparison theorem for such equations, but one can prove uniqueness directly by showing that any bounded solution of (4.9) coincides with the exponential utility indifference value process $C(B;\alpha)$. The main difficulty with (4.9) is that the presence of the compensated sum of jumps makes it very hard to derive any properties of the solution $Y$.

Note that both $\psi$ and $L$ in the BSDE (4.5) depend on the risk aversion parameter $\alpha$. We shall indicate this by writing $\psi(\alpha), L(\alpha)$.

**5. Dynamic and further properties of the indifference valuation.** In this section we derive further properties of the exponential utility indifference value process $C(B;\alpha)$. While some hold generally, others rely on the BSDE description in Theorem 13 and thus need continuity of $\mathbb{F}$. This will be specified if necessary so that the only standing assumptions in this section are that

$$S \text{ is locally bounded and } \mathbb{P}_{e,f} \neq \varnothing.$$

We first prove continuity of $C(B;\alpha)$ in $B$.

PROPOSITION 14. *Assume that $\mathbb{F}$ is continuous. If $(B^n)_{n\in\mathbb{N}}$ is a bounded sequence in $L^\infty$ such that $(B^n)$ converges to $B$ in probability for some $B \in L^\infty$, then for any $\gamma > 0$,*

$$(5.1) \quad \sup_{\alpha\in(0,\gamma]} \sup_{0\leq t\leq T} |C_t(B^n;\alpha) - C_t(B;\alpha)| \longrightarrow 0 \quad \text{in probability as } n\to\infty.$$

PROOF. We go back to the proof of Theorem 8 and work there with the pair $(S, Q^E)$ instead of $(M, R)$ and the two generators $(0, -\frac{\alpha}{2}, B^n)$ and $(0, -\frac{\alpha}{2}, B)$. The corresponding solutions are $(C(B^n; \alpha), \psi^n(\alpha), L^n(\alpha))$ and $(C(B; \alpha), \psi(\alpha), L(\alpha))$ by Theorem 13. From (3.11), we get

$$C_t(B^n; \alpha) - C_t(B; \alpha) = E_{Q^n(\alpha)}[B^n - B | \mathcal{F}_t],$$

where $Q^n(\alpha)$ is given by

$$dQ^n(\alpha) = \mathcal{E}\left(\frac{\alpha}{2}(L^n(\alpha) + L(\alpha))\right)_T dQ^E =: Z_T^n(\alpha) \, dQ^E.$$

The estimate (3.18) implies that

$$(5.2) \quad \sup_{\alpha \in (0, \gamma]} \sup_{n \in \mathbb{N}} \left\| \frac{\alpha}{2}(L^n(\alpha) + L(\alpha)) \right\|_{BMO(R)}^2 \leq \frac{\gamma}{2} \sup_{n \in \mathbb{N}} (e^{\|B^n\|_\infty} + e^{\|B\|_\infty})^2 < \infty,$$

and so there exists, by Theorem 3.1 of [24], an exponent $p \in (1, \infty)$ such that each $Z^n(\alpha)$ satisfies the reverse Hölder inequality $\mathcal{R}_p(Q^E)$, that is,

$$\sup_{0 \leq t \leq T} E_{Q^E}\left[\left(\frac{Z_T^n(\alpha)}{Z_t^n(\alpha)}\right)^p \Big| \mathcal{F}_t\right] \leq (c_p)^p$$

for a constant $c_p$. Note that because (5.2) is uniform in $n \in \mathbb{N}$ and $\alpha \in (0, \gamma]$, the same $p, c_p$ work for all these $n, \alpha$ simultaneously. Using now Bayes' rule and Hölder's inequality, we get

$$\sup_{\alpha \in (0, \gamma]} \sup_{0 \leq t \leq T} |C_t(B^n; \alpha) - C_t(B; \alpha)| = \sup_{\alpha \in (0, \gamma]} \sup_{0 \leq t \leq T} \left| E_{Q^E}\left[\frac{Z_T^n(\alpha)}{Z_t^n(\alpha)}(B^n - B) | \mathcal{F}_t\right] \right|$$

$$\leq c_p \sup_{0 \leq t \leq T} (E_{Q^E}[|B^n - B|^q | \mathcal{F}_t])^{1/q},$$

with $q \in (1, \infty)$ conjugate to $p$, and so (5.1) follows from Doob's maximal inequality. □

A natural assumption on a convex monetary utility functional $\Phi_t : L^\infty(\mathcal{F}_T) \to L^\infty(\mathcal{F}_t)$ is a continuity of the following form: If a bounded sequence $(B^n)_{n \in \mathbb{N}}$ in $L^\infty$ increases (or decreases), $P$-a.s. to some $B \in L^\infty$, then $\Phi_t(B^n)$ increases (or decreases), $P$-a.s. to $\Phi_t(B)$. This is one possible extension to the dynamic case of the semicontinuity requirements studied for static risk measures (or utility functionals); see, for instance, [14] or [10] for a recent conditional version. For the functional $\Phi_0 := C_0(\cdot; \alpha)$, the exponential utility indifference value at time 0, this continuity could be deduced from the recent work of Barrieu and El Karoui [2]; see their Theorem 3.6 and Proposition 5.3. However, Proposition 14 is stronger in that it provides such a result uniformly in $t \in [0, T]$ (and locally uniformly in $\alpha$ as well).

The next result holds generally, that is, without continuity of $\mathbb{F}$; see also Corollary 3.10 of [3].

PROPOSITION 15. *For each $\alpha \in (0,\infty)$, $C(B;\alpha)$ is time-consistent in the sense that, for any $B \in L^\infty$, we have*

(5.3)
$$C_\sigma(C_\tau(B;\alpha);\alpha) = C_\sigma(B;\alpha),$$
*P-a.s. for any stopping times $\sigma, \tau$ with $\sigma \leq \tau$.*

PROOF. Because $C_\tau(B';\alpha) = B'$ for any $\mathcal{F}_\tau$-measurable $B'$, we obtain from the dynamic programming equation (4.3) applied to $B' = C_\tau(B;\alpha)$ that

$$C_\sigma(B';\alpha) = \frac{1}{\alpha} \log \operatorname*{ess\,inf}_{\vartheta \in \Theta} E_{Q^E}[e^{\alpha(C_\tau(B;\alpha) - G_{\sigma,\tau}(\vartheta))} | \mathcal{F}_\sigma] = C_\sigma(B;\alpha), \qquad P\text{-a.s.} \quad \square$$

The financial interpretation of (5.3) is obvious: If we want to value the time $T$ payoff $B$ at time $\sigma$, we can either do this directly or first value $B$ at time $\tau \geq \sigma$ and then value the result $C_\tau(B;\alpha)$ at time $\sigma$. In both cases, the final valuation is the same. As emphasized by Musiela and Zariphopoulou [30], such a consistency property is highly desirable, and it is also known from the work of Rosazza Gianin [32] that a nice BSDE representation is usually sufficient to derive it. For more discussion and references on time-consistency aspects, we refer to [1].

As a direct consequence of Theorem 13 and Proposition 10, we also have the following:

PROPOSITION 16. *If $\mathbb{F}$ is continuous, the exponential utility indifference value $C_t(B;\alpha)$ is locally Lipschitz-continuous in $\alpha$, uniformly in $t$: For any $\gamma > 0$, we have*

$$\sup_{0 \leq t \leq T} |C_t(B;\alpha) - C_t(B;\alpha')| \leq K_\gamma |\alpha - \alpha'|, \qquad P\text{-a.s.}$$

*for all $\alpha, \alpha' \in (0, \gamma]$, where the constant $K_\gamma$ depends only on $\gamma$ and $B$.*

**6. Risk aversion asymptotics.** In this section we study the behavior of the exponential utility indifference value process as the risk aversion parameter $\alpha$ goes to 0 or $\infty$. Earlier results on some aspects of this have been obtained by Rouge and El Karoui [33], Becherer [3], Fujiwara and Miyahara [18] and Stricker [34], among others; see below for more detailed comments. As before, our standing assumptions in this section are that

$$S \text{ is locally bounded and } \mathbb{P}_{e,f} \neq \varnothing.$$

6.1. *Asymptotics for $\alpha \searrow 0$.* A simple adaptation of arguments from [34] gives the following:

THEOREM 17. *For each $B \in L^\infty$, we have*

(6.1) $\quad \lim_{\alpha \to 0} C_t(B; \alpha) = E_{Q^E}[B|\mathcal{F}_t] \qquad$ *uniformly in $t \in [0,T]$, $P$-a.s.*

*Moreover, we have the estimate*

(6.2) $\quad \sup_{0 \leq t \leq T} |C_t(B; \alpha) - E_{Q^E}[B|\mathcal{F}_t]| \leq \alpha\, const.(B), \qquad P\text{-a.s.}$

PROOF. With the notation $Z_{t,T} := Z_T/Z_t$, we know from Lemma 6 and the representations (1.7) and (1.5) that, for any $t \in [0,T]$, $\alpha \in (0,\infty)$ and $Q \in \mathbb{P}_{e,f}$,

$$E_{Q^E}[B|\mathcal{F}_t] \leq C_t(B;\alpha)$$
$$\leq E_Q[B|\mathcal{F}_t] - \frac{1}{\alpha}(E_Q[\log Z_{t,T}^Q|\mathcal{F}_t] - E_{Q^E}[\log Z_{t,T}^E|\mathcal{F}_t]).$$

Moreover, the representation (1.4) of $Z_T^E$ implies that

$$E_Q[\log Z_{t,T}^E|\mathcal{F}_t] = E_{Q^E}[\log Z_{t,T}^E|\mathcal{F}_t] \qquad \text{for any } Q \in \mathbb{P}_{e,f},$$

and we have

$$\log Z_{t,T}^Q - \log Z_{t,T}^E = \log(Z_{t,T}^Q/Z_{t,T}^E) = \log Z_{t,T}^{Q:Q^E},$$

where $Z^{Q:Q^E}$ denotes the density process of $Q$ with respect to $Q^E$. Bayes' rule and the Fenchel inequality $bz \leq \frac{1}{\alpha}(e^{\alpha b} + z \log z - 1)$ thus give

$$E_Q[B|\mathcal{F}_t] = E_{Q^E}[B Z_{t,T}^{Q:Q^E}|\mathcal{F}_t]$$
$$\leq \frac{1}{\alpha}(E_{Q^E}[e^{\alpha B}|\mathcal{F}_t] + E_{Q^E}[Z_{t,T}^{Q:Q^E} \log Z_{t,T}^{Q:Q^E}|\mathcal{F}_t] - 1)$$
$$= \frac{1}{\alpha}(E_{Q^E}[e^{\alpha B}|\mathcal{F}_t] + E_Q[\log Z_{t,T}^Q - \log Z_{t,T}^E|\mathcal{F}_t] - 1),$$

and so we get

$$\sup_{0 \leq t \leq T} |C_t(B;\alpha) - E_{Q^E}[B|\mathcal{F}_t]| \leq \sup_{0 \leq t \leq T} E_{Q^E}\left[\frac{e^{\alpha B} - 1}{\alpha} - B\Big|\mathcal{F}_t\right].$$

Because $B$ is bounded, we have $0 \leq \frac{e^{\alpha B}-1}{\alpha} - B \leq \frac{\alpha}{2}\|B\|_\infty^2 + \text{const.}\alpha^2$, $P$-a.s., and so (6.2) and (6.1) both follow. $\square$

REMARK. The convergence $\lim_{\alpha \to 0} C_t(B;\alpha) = E_{Q^E}[B|\mathcal{F}_t]$ has also been obtained by Rouge and El Karoui [33] for arbitrary (but fixed) $t$ in a Brownian filtration, and for $t = 0$ by Becherer [3] and Stricker [34] in a general setting and by Fujiwara and Miyahara [18] for geometric Lévy processes. Theorem 17 extends the argument by Stricker [34], who also gave the convergence rate of order $\alpha$, to provide a uniform result for all $t \in [0,T]$.

If $\mathbb{F}$ is continuous, an alternative proof of Theorem 17 goes via the BSDE description of $C(B;\alpha)$ in Theorem 13. In fact, taking conditional expectations between $t$ and $T$ in (4.5) and using (3.2) and the fact that $\int \psi(\alpha)\,dS$ and $L(\alpha)$ are $Q^E$-martingales yields

$$C_t(B;\alpha) = E_{Q^E}[B|\mathcal{F}_t] + \frac{\alpha}{2} E_{Q^E}[\langle L(\alpha)\rangle_T - \langle L(\alpha)\rangle_t|\mathcal{F}_t].$$

Hence, (6.2) follows from the estimate (3.15) in Lemma 9. We now prove that we also have convergence of the strategies $\psi(\alpha)$.

THEOREM 18. *Suppose that $\mathbb{F}$ is continuous and write the Galtchouk–Kunita–Watanabe decomposition of $B \in L^\infty$ under $Q^E$ as*

$$(6.3) \qquad V^E := E_{Q^E}[B|\mathbb{F}] = V_0^E + \int \psi^E\,dS + L^E.$$

*Then we have*

$$(6.4) \qquad \lim_{\alpha \to 0} \int \psi(\alpha)\,dS = \int \psi^E\,dS \quad \text{in } BMO(Q^E),$$

$$(6.5) \qquad \lim_{\alpha \to 0} L(\alpha) = L^E \quad \text{in } BMO(Q^E)$$

*and, more precisely, we even have*

$$(6.6) \qquad \left\| \int \psi(\alpha)\,dS - \int \psi^E\,dS \right\|_{BMO(Q^E)}^2 + \|L(\alpha) - L^E\|_{BMO(Q^E)}^2 \leq \alpha\,\mathrm{const.}(B).$$

PROOF. Since $\mathbb{F}$ is continuous, all processes below are continuous. Using (6.3) and Theorem 13, we obtain from Itô's formula, omitting the arguments $B$ and $\alpha$ for the moment, that

$$(6.7) \quad \begin{aligned}(C_T - V_T^E)^2 &= (C_t - V_t^E)^2 - 2\int_t^T (C_s - V_s^E)\frac{\alpha}{2}\,d\langle L\rangle_s \\ &\quad + \int_t^T (\psi_u - \psi_u^E)^{\mathrm{tr}}\,d\langle S\rangle_u(\psi_u - \psi_u^E) \\ &\quad + \int_t^T d\langle L - L^E\rangle_s \\ &\quad + 2\int_t^T (C_u - V_u^E)\,d\left(\int(\psi - \psi^E)\,dS + L - L^E\right)_u.\end{aligned}$$

Since $V^E$ is a bounded $Q^E$-martingale, $\int \psi^E \, dS$ and $L^E$ are in $BMO(Q^E)$ and thus $Q^E$-martingales. Hence, the last term in (6.7) is like its integrator a $Q^E$-martingale because the integrand is bounded. Taking conditional expectations and using $C_T(B) = B = V_T^E$ yields

$$E_{Q^E}\left[\int_t^T (\psi_u(\alpha) - \psi_u^E)^{\mathrm{tr}} \, d\langle S\rangle_u (\psi_u(\alpha) - \psi_u^E) \Big| \mathcal{F}_t\right]$$

$$+ E_{Q^E}\left[\int_t^T d\langle L(\alpha) - L^E\rangle_s \Big| \mathcal{F}_t\right] + (C_t(B;\alpha) - V_t^E)^2$$

$$= \alpha E_{Q^E}\left[\int_t^T (C_s(B;\alpha) - V_s^E) \, d\langle L(\alpha)\rangle_s \Big| \mathcal{F}_t\right]$$

$$\leq 2\|B\|_\infty \alpha E_{Q^E}[\langle L(\alpha)\rangle_T - \langle L(\alpha)\rangle_t | \mathcal{F}_t]$$

$$\leq 2\|B\|_\infty \alpha \sup_{\alpha \in (0,\infty)} \|L(\alpha)\|_{BMO(Q^E)}^2 \qquad \text{uniformly in } t.$$

Hence, (6.4)–(6.6) all follow from (3.16), and we also again recover (6.1). □

Loosely speaking, the interpretation of Theorem 18 is that, in the small risk aversion limit, exponential indifference hedging converges to risk-minimization under the minimal entropy martingale measure $Q^E$. To see this, note that the integrand $\psi^E$ in the decomposition (6.3) of $B$ is (the risky asset component of) the strategy which is risk-minimizing in the sense of Föllmer and Sondermann [15] with respect to $Q^E$. Hence, Theorem 18 says that, for vanishing risk aversion $\alpha$, the gains process $\int \psi(\alpha) \, dS$ from the $\alpha$-optimal strategy for exponential utility indifference valuation converges to the gains process from the $Q^E$-risk-minimizing strategy. As in Theorem 17, we even obtain a convergence rate.

REMARK. The convergence in (6.4) was conjectured by D. Becherer in private discussions with one of the authors. Theorem 18 also explains the observation made after Corollary 4.3 of Young [35] that, in a particular model for valuing catastrophe bonds by exponential utility indifference, formally setting $\alpha = 0$ reproduces an earlier alternative approach.

6.2. *Asymptotics for $\alpha \nearrow \infty$.* Our last contribution is a study of the large risk aversion asymptotics of $C(B;\alpha)$. To that end, we recall the *super-replication price process*

$$C_t^*(B) := \operatorname*{ess\,sup}_{Q \in \mathbb{P}_e} E_Q[B|\mathcal{F}_t], \qquad 0 \leq t \leq T,$$

where we can and do choose an RCLL version. By the optional decomposition theorem (see [12] or [27]), $C^*(B)$ is the smallest RCLL process with final value $B$ at time $T$ which is a $Q$-supermartingale for all $Q \in \mathbb{P}_e$, and it admits a decomposition

$$(6.8) \qquad C^*(B) = C_0^*(B) + \int \psi^* \, dS - K^*,$$

where $\psi^*$ is an $\mathbb{R}^d$-valued predictable $S$-integrable process and $K^*$ is an optional increasing process null at 0. In general, $K^*$ is neither unique nor continuous; see Example 1 of [12]. But if the filtration is continuous, $K^*$ is actually predictable, hence, unique by the Doob–Meyer decomposition theorem, and because $C^*(B)$ is bounded, that result then also implies that $K^*$ is $Q^E$-integrable and $\psi^*$ is in $\Theta$.

From part (2) of Lemma 6, we know that

$$(6.9) \qquad C_t(B;\alpha) \leq C_t^*(B), \qquad P\text{-a.s. for each } t \in [0,T].$$

Moreover, we also have

$$(6.10) \qquad C_t^*(B) = \operatorname*{ess\,sup}_{Q \in \mathbb{P}_{e,f}} E_Q[B|\mathcal{F}_t], \qquad P\text{-a.s. for each } t \in [0,T].$$

In fact, Bayes' rule gives

$$E_Q[B|\mathcal{F}_t] = E_P[Z_T^Q B|\mathcal{F}_t]/E_P[Z_T^Q|\mathcal{F}_t] \qquad \text{for } Q \in \mathbb{P}_e,$$

and by Theorem 1.1 and Corollary 1.3 of [22], the set $\{Z_T^Q|Q \in \mathbb{P}_e\} \cap L^\infty(P) \subseteq \{Z_T^Q|Q \in \mathbb{P}_{e,f}\}$ is dense in $\{Z_T^Q|Q \in \mathbb{P}_e\}$ for the $L^1(P)$-norm. Since $B \in L^\infty(P)$, (6.10) readily follows.

For the next result, we need some notation. Let $D = (D_t)_{0 \leq t \leq T}$ be an increasing predictable RCLL process null at 0 such that $\langle S^i, S^k \rangle \ll D$ for all $i,k = 1, \ldots, d$. We choose $D$ strictly increasing and bounded (uniformly in $t,\omega$); for instance, $D_t := t + \tanh(\sum_{i=1}^d \langle S^i \rangle_t)$ will do. If $S$ is continuous, we can and do choose $D$ continuous as well. We define the $(d \times d)$ matrix-valued predictable process $\underline{\Sigma} = (\underline{\Sigma}_t)_{0 \leq t \leq T}$ by $d\langle S \rangle_t = \underline{\Sigma}_t \, dD_t$ and the finite measure $\mu$ on $\overline{\Omega} := \Omega \times [0,T]$ by $\mu := Q^E \otimes D$. Then we have

$$\left\langle \int \vartheta \, dS, \int \vartheta' \, dS \right\rangle = \int \vartheta^{\mathrm{tr}} \underline{\Sigma} \vartheta' \, dD \qquad \text{for } \vartheta, \vartheta' \in L(S)$$

and

$$(6.11) \qquad \left\| \int_0^T \vartheta_u \, dS_u \right\|_{L^2(Q^E)}^2 = \left\| \left\langle \int \vartheta \, dS \right\rangle_T \right\|_{L^1(Q^E)} = \|\vartheta^{\mathrm{tr}} \underline{\Sigma} \vartheta\|_{L^1(\mu)}$$

if $\int \vartheta \, dS$ is square-integrable under $Q^E$. For $d = 1$, we do not need all this notation since we can take $D = \langle S \rangle$ and $\underline{\Sigma} \equiv 1$; the measure $\mu$ is then the Doléans measure of $\langle S \rangle$ under $Q^E$.

THEOREM 19. *Assume that $\mathbb{F}$ is continuous. Fix $B \in L^\infty$ and any stopping time $\tau$. Then*:

(1) $P[\lim_{\alpha \to \infty} C_t(B; \alpha) = C_t^*(B)$ for all $t \in [0, T]] = 1$. (*This is true even without continuity of $\mathbb{F}$ or $S$.*)
(2) $\lim_{\alpha \to \infty} C_\tau(B; \alpha) = C_\tau^*(B)$ *strongly in* $L^r(Q^E)$ *for every* $r \in [1, \infty)$.
(3) $\lim_{\alpha \to \infty} \int_0^\tau \psi_u(\alpha) \, dS_u = \int_0^\tau \psi_u^* \, dS_u$ *weakly in* $L^r(Q^E)$ *for every* $r \in [1, \infty)$.
(4) $\lim_{\alpha \to \infty} \frac{\alpha}{2} \langle L(\alpha) \rangle_\tau = K_\tau^*$ *weakly in* $L^r(Q^E)$ *for every* $r \in [1, \infty)$.
(5) $\lim_{\alpha \to \infty} ((\psi(\alpha) - \psi^*)^{\text{tr}} \underline{\Sigma} (\psi(\alpha) - \psi^*))^{1/2} = 0$ *strongly in* $L^p(\mu)$ *for every* $p \in [1, 2)$.

PROOF. (a) The first part of the argument is almost as in [33]. From (6.9), (1.7) and (1.5), we have, for any $Q \in \mathbb{P}_{e,f}$, that

$$C_t^*(B) \geq C_t(B; \alpha) \geq E_Q[B|\mathcal{F}_t] - \frac{1}{\alpha}\left(E_Q\left[\log\frac{Z_T^Q}{Z_t^Q}\bigg|\mathcal{F}_t\right] - E_{Q^E}\left[\log\frac{Z_T^E}{Z_t^E}\bigg|\mathcal{F}_t\right]\right).$$

Letting $\alpha \to \infty$ and using (6.10) yields $\lim_{\alpha \to \infty} C_t(B; \alpha) = C_t^*(B)$, $P$-a.s. for each $t \in [0, T]$. Then (1) follows because $C(B; \alpha)$ and $C^*(B)$ are both right-continuous, and (2) then follows because all these processes are uniformly bounded by $\|B\|_\infty$. Clearly, this argument does not use the continuity of $\mathbb{F}$ or $S$.

(b) We already know that $C^*(B)$ and each $C(B; \alpha)$ are RCLL $Q^E$-supermartingales; see the remark following Proposition 12. Because we also have the convergence in (2) and a uniform bound $\|B\|_\infty$ on all these processes, Theorem VII.18 of [9] implies that at each stopping time, the $Q^E$-compensators converge weakly in $L^1(Q^E)$ as $\alpha \to \infty$. This still does not need continuity of $\mathbb{F}$ or $S$, but it also does not lead us very far because we cannot identify the compensators in general.

(c) Now assume that $\mathbb{F}$ is continuous. Then $C(B; \alpha)$ can be written as

$$(6.12) \qquad C(B; \alpha) = C_0(B; \alpha) + \int \psi(\alpha) \, dS + L(\alpha) - \frac{\alpha}{2}\langle L(\alpha) \rangle$$

by Theorem 13. From (6.12) and (6.8), we can therefore identify the $Q^E$-compensators as $\frac{\alpha}{2}\langle L(\alpha) \rangle$ and $K^*$, respectively, so that (b) gives

$$(6.13) \qquad \lim_{\alpha \to \infty} \frac{\alpha}{2} \langle L(\alpha) \rangle_\tau = K_\tau^* \qquad \text{weakly in } L^1(Q^E).$$

Due to (3.17) in Lemma 9, $L(\alpha)$ converges to 0 in $BMO(Q^E)$ as $\alpha \to \infty$, and this implies

$$(6.14) \qquad \lim_{\alpha \to \infty} L_\tau(\alpha) = 0 \qquad \text{strongly in } L^2(Q^E).$$

By combining (6.14) and (6.13) with (2) and (6.12) and (6.8), we obtain

$$(6.15) \qquad \lim_{\alpha \to \infty} \int_0^\tau \psi_u(\alpha) \, dS_u = \int_0^\tau \psi_u^* \, dS_u \qquad \text{weakly in } L^1(Q^E).$$

Hence, (3) follows from (6.15) if the family $\{\int_0^\tau \psi_u(\alpha)\,dS_u | \alpha \in (0,\infty)\}$ is bounded in $L^r(Q^E)$ for every $r \in [1,\infty)$. In view of (2) and (6.12), each of the families $\{\frac{\alpha}{2}\langle L(\alpha)\rangle_\tau | \alpha \in (0,\infty)\}$ and $\{L_\tau(\alpha) | \alpha \in (0,\infty)\}$ is then bounded in $L^r(Q^E)$ if and only if the other one is, and so (4) follows from (6.13) if $\{L_\tau(\alpha) | \alpha \in (0,\infty)\}$ is bounded in $L^r(Q^E)$ for every $r \in [1,\infty)$.

(d) For $N^\alpha \in \{\int \psi(\alpha)\,dS, L(\alpha)\}$, the energy inequalities give, for each $n \in \mathbb{N}$,

$$\sup_{0 \le t \le T} E_{Q^E}[(\langle N^\alpha\rangle_T - \langle N^\alpha\rangle_t)^n | \mathcal{F}_t] \le n! \|N^\alpha\|^{2n}_{BMO(Q^E)};$$

see [24], page 28. Using the Burkholder–Davis–Gundy inequalities and the estimates (3.14) and (3.16) in Lemma 9, applied with $(M,R) = (S, Q^E)$, thus yields

$$\sup_{\alpha \in (0,\infty)} E_{Q^E}\left[\left(\sup_{0 \le t \le T} |N^\alpha_t|\right)^{2n}\right] \le \sup_{\alpha \in (0,\infty)} \text{const.}(n) E_{Q^E}[(\langle N^\alpha\rangle_T)^n]$$

$$\le \text{const.}(n) n! \left(\sup_{\alpha \in (0,\infty)} \|N^\alpha\|_{BMO(Q^E)}\right)^{2n} < \infty.$$

Hence, $\{N^\alpha_\tau | \alpha \in (0,\infty)\}$ is bounded in $L^r(Q^E)$ for every $r \in [1,\infty)$, as desired in (c).

(e) To prove (5), we set $\eta(\alpha) := \psi(\alpha) - \psi^*$ and note from (3) that $\{\int_0^\tau \eta_u(\alpha)\,dS_u | \alpha \in (0,\infty)\}$ is bounded in $L^r(Q^E)$ for every $r \in [1,\infty)$. In view of (6.11), this means, for $r = 2$, that

$$\sup_{\alpha \in (0,\infty)} \|\eta(\alpha)^{\text{tr}} \underline{\Sigma} \eta(\alpha)\|_{L^1(\mu)} = \sup_{\alpha \in (0,\infty)} \left\|\int_0^T \eta_u(\alpha)\,dS_u\right\|^2_{L^2(Q^E)} < \infty$$

so that the family $\{(\eta(\alpha)^{\text{tr}} \underline{\Sigma}\eta(\alpha))^{1/2} | \alpha \in (0,\infty)\}$ is bounded in $L^2(\mu)$. Hence, (5) follows as soon as we prove that

(6.16) $$\lim_{\alpha \to \infty} \eta(\alpha)^{\text{tr}} \underline{\Sigma} \eta(\alpha) = 0 \quad \text{in } \mu\text{-measure}.$$

(f) The proof of (6.16) is a slight variation of an argument due to Peng [31]. We first apply Itô's formula and use (6.12) and (6.8), suppressing for the moment all arguments $\alpha$ and $B$, to obtain, for any stopping times $\sigma \le \rho$,

(6.17)
$$\begin{aligned}(C^*_\rho - C_\rho)^2 &= (C^*_\sigma - C_\sigma)^2 + 2\int_\sigma^\rho (C^*_{s-} - C_{s-})\left(\frac{\alpha}{2}\,d\langle L\rangle_s - dK^*_s\right) \\ &\quad + \int_\sigma^\rho (\psi^*_u - \psi_u)^{\text{tr}}\,d\langle S\rangle_u (\psi^*_u - \psi_u) + \langle L\rangle_\rho - \langle L\rangle_\sigma + [K^*]_\rho - [K^*]_\sigma \\ &\quad + 2\int_\sigma^\rho (C^*_{u-} - C_{u-})\,d\left(\int(\psi^* - \psi)\,dS - L\right)_u.\end{aligned}$$

The last term is a $Q^E$-martingale because the integrand is bounded and the integrator is a $Q^E$-martingale due to $\psi^* \in \Theta$ and Lemma 9. Moreover, $C = C(B;\alpha)$ is continuous by Theorem 13, and (6.8) gives $\Delta K^* = -\Delta C^*$ because $S$ is continuous. Hence, $[K^*] = \sum (\Delta K^*_s)^2 = -\int \Delta C^* \, dK^*$ and, therefore,

$$\int (C^*_- - C_-)\left(\frac{\alpha}{2} d\langle L \rangle - dK^*\right) + [K^*] = \int (C^* - C)\left(\frac{\alpha}{2} d\langle L \rangle - dK^*\right).$$

Adding and subtracting $[K^*]_\rho - [K^*]_\sigma$ in (6.17) and taking expectations therefore yields

$$E_{Q^E}\left[\int_\sigma^\rho (\psi^*_u - \psi_u)^{\mathrm{tr}} \, d\langle S\rangle_u \, (\psi^*_u - \psi_u)\right]$$
$$+ E_{Q^E}[\langle L\rangle_\rho - \langle L\rangle_\sigma] + E_{Q^E}[(C^*_\sigma - C_\sigma)^2]$$
$$= E_{Q^E}[(C^*_\rho - C_\rho)^2] + 2E_{Q^E}\left[\int_\sigma^\rho (C^*_s - C_s)\left(dK^*_s - \frac{\alpha}{2} d\langle L\rangle_s\right)\right]$$
$$+ E_{Q^E}[[K^*]_\rho - [K^*]_\sigma]$$
$$\leq E_{Q^E}[(C^*_\rho - C_\rho)^2] + 2E_{Q^E}\left[\int_0^T (C^*_s - C_s) \, dK^*_s\right] + E_{Q^E}[[K^*]_\rho - [K^*]_\sigma]$$

because $C^* - C \geq 0$ by (6.9). On the left-hand side, the middle term goes to 0 as $\alpha \to \infty$ by (3.17), and the last term goes to 0 as well, due to (1). On the right-hand side, the first term goes to 0 as $\alpha \to \infty$ due to (1) and the second by using (1) and dominated convergence, because $K^*_T \in L^1(Q^E)$. Since $\eta(\alpha) = \psi(\alpha) - \psi^*$, we thus obtain that

$$(6.18) \quad \begin{aligned} \limsup_{\alpha \to \infty} E_{Q^E}\left[\int_\sigma^\rho \eta(\alpha)^{\mathrm{tr}} \underline{\Sigma} \eta(\alpha) \, dD\right] \\ = \limsup_{\alpha \to \infty} E_{Q^E}\left[\int_\sigma^\rho \eta_u(\alpha)^{\mathrm{tr}} \, d\langle S\rangle_u \, \eta_u(\alpha)\right] \\ \leq E_{Q^E}[[K^*]_\rho - [K^*]_\sigma] \end{aligned}$$

for all stopping times $\sigma \leq \rho$.

Now we use Lemma 20 below (with $A = K^*$ and $\beta = D$) to obtain, for any $\delta, \varepsilon > 0$, finitely many pairwise disjoint intervals $]]\sigma_k, \tau_k]]$, $k = 0, 1, \ldots, N$, such that $0 < \sigma_k \leq \tau_k \leq T$ and

$$(6.19) \quad \mu\left(\bigcup_{k=0}^N ]]\sigma_k, \tau_k]]\right) = E_{Q^E}\left[\sum_{k=0}^N (D_{\tau_k} - D_{\sigma_k})\right]$$
$$\geq E_{Q^E}[D_T] - \frac{\varepsilon}{2} = \mu(\overline{\Omega}) - \frac{\varepsilon}{2},$$

$$(6.20) \quad \sum_{k=0}^N E_{Q^E}\left[\sum_{\sigma_k < t \leq \tau_k} (\Delta K^*_t)^2\right] \leq \frac{\delta\varepsilon}{2}.$$

Note that $E_{Q^E}[(K_T^*)^2] < \infty$ follows from (4). Applying the estimate (6.18) for each $\sigma = \sigma_k$, $\rho = \tau_k$ and taking the sum from $k = 0$ to $N$, we have from (6.20) that

$$\limsup_{\alpha \to \infty} \sum_{k=0}^{N} E_{Q^E}\left[\int_{\sigma_k}^{\tau_k} \eta(\alpha)^{\mathrm{tr}} \underline{\Sigma} \eta(\alpha)\, dD\right] \leq \sum_{k=0}^{N} E_{Q^E}\left[\sum_{\sigma_k < t \leq \tau_k} (\Delta K_t^*)^2\right] \leq \frac{\delta\varepsilon}{2}.$$

Thus, there exists some $\alpha_0(\delta, \varepsilon)$ such that, for all $\alpha \geq \alpha_0(\delta, \varepsilon)$, we have

$$\sum_{k=0}^{N} E_{Q^E}\left[\int_{\sigma_k}^{\tau_k} \eta(\alpha)^{\mathrm{tr}} \underline{\Sigma} \eta(\alpha)\, dD\right] \leq \frac{\delta\varepsilon}{2},$$

which implies by Markov's inequality that

$$\mu\left(\left(\bigcup_{k=0}^{N} ]]\sigma_k, \tau_k]]\right) \cap \{\eta(\alpha)^{\mathrm{tr}} \underline{\Sigma} \eta(\alpha) \geq \delta\}\right) \leq \frac{\varepsilon}{2}.$$

Combining this with (6.19) implies that

$$\mu(\{\eta(\alpha)^{\mathrm{tr}} \underline{\Sigma} \eta(\alpha) \geq \delta\}) \leq \varepsilon \qquad \text{for all } \alpha \geq \alpha_0(\delta, \varepsilon)$$

so that $\eta(\alpha)^{\mathrm{tr}} \underline{\Sigma} \eta(\alpha)$ converges to 0 in $\mu$-measure. This completes the proof. $\square$

REMARKS. (1) The pointwise convergence in (1) of Theorem 19 has also been given by Rouge and El Karoui [33], although it is not quite clear from their proof how (6.10) comes in. In addition to a uniform result in $t$, we also provide here in (3) and (5) the convergence of the strategies and in (4) of the residual terms in the BSDE for $C(B; \alpha)$.

(2) To the best of our knowledge, Theorem 19 is the first result in continuous time on the convergence of strategies in utility indifference valuation. For related work in a one-period model, see [6].

In the proof of Theorem 19, we have used the following technical result originally due to Peng [31] for the case $\beta_t = t$.

LEMMA 20. *Suppose that the filtration $\mathbb{F}$ is continuous. Let $A = (A_t)_{0 \leq t \leq T}$ be an increasing RCLL process with $A_0 = 0$ and $E[A_T^2] < \infty$, and let $\beta = (\beta_t)_{0 \leq t \leq T}$ be a continuous increasing process with $\beta_0 = 0$ and $E[\beta_T] < \infty$. Then for any $\delta, \varepsilon > 0$, there exist finitely many stopping times $\sigma_k, \tau_k$, $k = 0, 1, \ldots, N$, with $0 < \sigma_k \leq \tau_k \leq T$ and such that*

(i) $\qquad ]]\sigma_i, \tau_i]] \cap ]]\sigma_k, \tau_k]] = \varnothing \qquad \text{for } i \neq k,$

(ii) $\qquad E\left[\sum_{k=0}^{N} (\beta_{\tau_k} - \beta_{\sigma_k})\right] \geq E[\beta_T] - \varepsilon,$

(iii) $$\sum_{k=0}^{N} E\left[\sum_{\sigma_k < t \leq \tau_k} (\Delta A_t)^2\right] \leq \delta.$$

PROOF. This is done almost exactly as in [31]. Continuity of $\mathbb{F}$ ensures that all stopping times are predictable, hence, foretellable, so that Lemma A.2 of [31] still holds. Continuity of $\beta$ guarantees that we can obtain (ii) as in [31]. □

**Acknowledgments.** M. Schweizer thanks Dirk Becherer and Susanne Klöppel for helpful discussions, and Christian Bender for suggesting to use the result by Peng [31].

DEPARTMENT OF PROBABILITY THEORY
  AND MATHEMATICAL STATISTICS
A. RAZMADZE MATHEMATICAL INSTITUTE
1 ALEXIDZE ST.
TBILISI 0193
GEORGIA
E-MAIL: mania@rmi.acnet.ge

DEPARTEMENT MATHEMATIK
ETH ZÜRICH
ETH-ZENTRUM, HG G28.2
CH-8092 ZÜRICH
SWITZERLAND
E-MAIL: martin.schweizer@math.ethz.ch